\documentclass[12pt]{amsart}
\usepackage{latexsym, amssymb, amsfonts, amscd}

\theoremstyle{plain}
\newtheorem{Thm}{Theorem}[section]

\newtheorem{Main}[Thm]{Main Theorem}

\newtheorem{Prop}[Thm]{Proposition}
\newtheorem{Lem}[Thm]{Lemma}
\newtheorem{Cl}[Thm]{Claim}

\theoremstyle{definition}
\newtheorem{Def}[Thm]{Definition}
\newtheorem{Rem}[Thm]{Remark}

\newtheorem{Emp}[Thm]{}

\newtheorem{Not}[Thm]{Notation}

\numberwithin{equation}{section}

\newcommand{\ov}{\overline}

\newcommand{\fq}{\B{F}_q}

\newcommand{\B}[1]{\mathbb#1}
\newcommand{\cal}[1]{\mathcal{#1}}
\newcommand{\C}[1]{\cal#1}

\newcommand{\sr}{\operatorname{sr}}
\newcommand{\ssc}{\operatorname{sc}}
\newcommand{\der}{\operatorname{der}}
\newcommand{\nr}{\operatorname{nr}}
\newcommand{\tu}{\operatorname{tu}}
\newcommand{\tn}{\operatorname{tn}}

\newcommand{\isom}{\overset {\thicksim}{\to}}

\newcommand{\Om}{\Omega}

\newcommand{\si}{\sigma}
\newcommand{\Si}{\Sigma}

\newcommand{\lra}{\longrightarrow}

\newcommand{\hra}{\hookrightarrow}
\newcommand{\wt}{\widetilde}
\newcommand{\wh}{\widehat}
\newcommand{\Gm}{\Gamma}

\newcommand{\gm}{\gamma}
\newcommand{\ka}{\kappa}
\newcommand{\dt}{\delta}
\newcommand{\Dt}{\Delta}

\newcommand{\m}{^{\times}}

\newcommand{\un}{\operatorname{un}}
\newcommand{\nil}{\operatorname{nil}}

\newcommand{\la}{\lambda}

\newcommand{\rl}[1]{Lemma \ref{L:#1}}
\newcommand{\rn}[1]{Notation \ref{N:#1}}
\newcommand{\rcl}[1]{Claim \ref{C:#1}}
\newcommand{\rp}[1]{Proposition \ref{P:#1}}
\newcommand{\rr}[1]{Remark \ref{R:#1}}

\newcommand{\re}[1]{\ref{E:#1}}

\newcommand{\rt}[1] {Theorem \ref{T:#1}}
\newcommand{\rd}[1]{Definition \ref{D:#1}}

\newcommand{\inv}{\operatorname{inv}}
\newcommand{\pr}{\operatorname{pr}}
\newcommand{\Ker}{\operatorname{Ker}}
\newcommand{\Out}{\operatorname{Out}}
\newcommand{\Coker}{\operatorname{Coker}}
\newcommand{\Norm}{\operatorname{Norm}}

\newcommand{\Aut}{\operatorname{Aut}}

\newcommand{\Ind}{\operatorname{Ind}}

\newcommand{\Ad}{\operatorname{Ad}}

\newcommand{\Tr}{\operatorname{Tr}}

\newcommand{\ad}{\operatorname{ad}}

\newcommand{\End}{\operatorname{End}}

\newcommand{\rk}{\operatorname{rk}}

\newcommand{\lan}{\langle}
\newcommand{\ran}{\rangle}

\begin{document}

\title[Endoscopic decomposition]%
{Endoscopic decomposition of characters of certain cuspidal representations}

\author{David Kazhdan}
\author{Yakov Varshavsky}
\address{Institute of Mathematics\\
Hebrew University\\
Givat-Ram, Jerusalem,  91904\\
Israel}
\email{kazhdan@math.huji.ac.il, vyakov@math.huji.ac.il }

\thanks{The work of the second author was supported by 
THE ISRAEL SCIENCE FOUNDATION (Grant No. 38/01-1)}
\date{\today}

\begin{abstract}
We construct an endoscopic decomposition for local $L$-packets 
associated to irreducible cuspidal Deligne--Lusztig 
representations. Moreover, the obtained decomposition is compatible with inner twistings.
\end{abstract}
\maketitle

\section{Introduction}

Let $E$ be a local non-archimedean field, with ring of integers $\C{O}$ and residue field 
$\fq$ of characteristic $p$. We denote by $\Gm\supset W\supset I$ the absolute Galois, 
the Weil and the inertia groups of $E$. Let $G$ be a reductive group over 
$E$, ${}^LG=\wh{G}\rtimes W$ its complex Langlands dual group, and $\C{D}(G(E))$ the 
space of invariant distributions on $G(E)$. 

Every admissible homomorphism $\la:W\to {}^LG$ (see \cite[$\S$ 10]{Ko1}) 
gives rise to a finite group $S_{\la}:=\pi_0(Z_{\wh{G}}(\la)/Z(\wh{G})^{\Gm})$,
where $Z_{\wh{G}}(\la)$ is the centralizer of $\la(W)$ in $\wh{G}$.
Every conjugacy class $\ka$ of $S_{\la}$ defines an endoscopic
subspace $\C{D}_{\ka,\la}(G(E))\subset\C{D}(G(E))$. For simplicity, we will restrict 
ourselves to the elliptic case, where $\la(W)$ does not lie in any proper Levi subgroup of ${}^LG$.

Langlands conjectured that every elliptic $\la$ corresponds to a finite set $\Pi_{\la}$,  
called an $L$-packet, of cuspidal irreducible representations of $G(E)$. 
Moreover, the subspace $\C{D}_{\la}(G(E))\subset\C{D}(G(E))$, generated by characters 
$\{\chi(\pi)\}_{\pi\in\Pi_{\la}}$, should have an endoscopic decomposition. 
More precisely, it is expected (\cite[IV, 2]{La})
that there exists a basis $\{a_{\pi}\}_{\pi\in \Pi_{\la}}$ of the space of central functions on  
$S_{\la}$ such that 
$\chi_{\ka,\la}:=\sum_{\pi\in \Pi_{\la}}a_{\pi}(\ka) \chi(\pi)$ belongs to 
$\C{D}_{\ka,\la}(G)$ for every conjugacy class $\ka$ of $S_{\la}$. 

The goal of this paper is to construct the endoscopic decomposition of $\C{D}_{\la}(G(E))$ for
tamely ramified $\la$'s such that $Z_{\wh{G}}(\la(I))$ is a maximal torus. In this case, $G$ splits over
an unramified extension of $E$, and $\la$ factors through ${}^LT\hra{}^LG$ for an elliptic 
unramified maximal torus $T$ of $G$.
By the local Langlands correspondence for tori (\cite{La2}), a homomorphism 
$\la:W\to {}^LT$ defines a tamely ramified homomorphism $\theta:T(E)\to\B{C}\m$. 
Each $\ka\in S_{\la}=\wh{T}^{\Gm}/Z(\wh{G})^{\Gm}$ gives rise to an elliptic endoscopic 
datum $\C{E}_{\ka,\la}$ of $G$, while characters of $S_{\la}$ are in bijection with 
conjugacy classes of embeddings $T\hra G$, stably conjugate to the inclusion. 
Therefore each character $a$ of $S_{\la}$ gives rise to an irreducible cuspidal 
representation 
$\pi_{a,\la}$ of $G(E)$ (denoted by $\pi_{a,\theta}$ in \rn{repr}). 

%

Our main result asserts, for fields $E$ of sufficiently large residual characteristic, 
that each $\chi_{\ka,\la}:=\sum_{a}a(\ka) \chi(\pi_{a,\la})$ is 
$\C{E}_{\ka,\la}$-stable. Moreover, the resulting endoscopic decomposition of 
$\C{D}_{\la}(G(E))$ is compatible with inner twistings. 
For simplicity, we restrict ourselves to local fields of characteristic zero,  
while the case of positive characteristic follows by approximation (see \cite{Ka3,De}). 


Our argument goes as follows. First we prove the stability 
of the restriction of $\chi_{\ka,\la}$ to the subset of topologically unipotent elements of 
$G(E)$. If $p$ is sufficiently large, this assertion reduces to 
the analogous assertion about distributions on the Lie algebra. 
Now the stability follows from a combination of a Springer hypothesis \cite{Ka1}
and a generalization of a theorem of Waldspurger \cite{Wa}.
To prove the result in general, we use the topological Jordan decomposition (\cite{Ka2}).


When this work was in the process of writing, we have heard that S. DeBacker and  M. Reeder
obtained similar results.

\vskip 4truept
\centerline{\bf Notation and Conventions.}
\vskip 4truept

In addition to the notation introduced above, we use the following conventions:

  For a reductive group $G$, always assumed to be connected, we denote by $Z(G)$, $G^{\ad}$, $G^{\der}$, 
$G^{\ssc}$, $G_{\dt}$
and $G^{\sr}$ the center of $G$, the adjoint group of $G$, the derived group of $G$, 
the simply connected covering of $G^{\der}$, the centralizer of $\dt\in G$, and 
the set of strongly regular semisimple elements of  $G$ (that is, the set of $\dt\in G$ such that 
$G_{\dt}\subset G$ is a maximal torus) respectively.

  Denote by $\C{G}, \C{T},$ and $\C{L}$ Lie algebras of the algebraic groups $G, T$ and
$L$.

  Let $E$ be a local non-archimedean field of characteristic zero,  $\ov{E}$
a fixed algebraic closure of $E$, and $E^{\nr}$ a maximal unramified extension of $E$ in $\ov{E}$.

  For a reductive group $G$ (resp. its Lie algebra $\C{G}$) over $E$, we denote
by  $\C{S}(G(E))$ (resp.  $\C{S}(\C{G}(E))$) the space of locally constant measures with compact support.
We denote by $\C{D}(G(E))$ (resp. $\C{D}(\C{G}(E))$) the space of  invariant distributions on $G(E)$ 
(resp. $\C{G}(E)$), namely $G(E)$-invariant linear functionals on
$\C{S}(G(E))$ (resp  $\C{S}(\C{G}(E))$), where  $G(E)$ acts by conjugation.
Whenever necessary we equip $G(E)$ and $\C{G}(E)$ with invariant measures, denoted by $\mu$, defined
by a translation-invariant top degree differential form on $G$.
We denote by $G(E)_{\tu}$ (resp. $\C{G}(E)_{\tn}$) the set of topologically unipotent 
(resp. topologically nilpotent) elements of $G(E)$ (resp. $\C{G}(E)$).
Finally, we denote by $\rk(G)$ the rank of $G$ over $E$, and put $e(G):=(-1)^{\rk(G^{\ad})}$.
Note that our sign $e(G)$ differs from that defined by Kottwitz. 

\section{Formulation of the Main result}

\begin{Emp} \label{E:finite}
Let $L$ be a connected reductive group over $\fq$, and
$\ov{a}:\ov{T}\hra L$ an embedding of a maximal elliptic torus of $L$.
Following Deligne and Lusztig,  we associate to every   
character $\ov{\theta}:\ov{T}(\fq)\to\B{C}\m$ in general position an irreducible cuspidal representation 
$\rho_{\ov{a},\ov{\theta}}$ of $L(\fq)$ (see \cite[Prop. 7.4 and Thm. 8.3]{DL}).
\end{Emp}

\begin{Emp}
There is an equivalence of categories $T\mapsto \ov{T}$ between tori over $E$ splitting over 
$E^{\nr}$ and tori over $\fq$. Every such $T$ has a canonical $\C{O}$-structure.
\end{Emp}

\begin{Not} \label{N:repr}

a) Let $G$ be a reductive group over $E$, $T$ a torus over $E$ splitting over $E^{\nr}$, and  
$a:T\hra G$ an embedding of a maximal elliptic torus of $G$. Then 
$G$ splits over $E^{\nr}$, and $a(T(\C{O}))$ lies in a unique parahoric subgroup $G_{a}$ of $G(E)$. 
Let $G_{a^+}$ be the pro-unipotent radical of $G_a$.
Then there exists a canonical reductive group $L_a$ over $\fq$ 
with an identification $L_a(\fq)=G_a/G_{a^+}$.
Moreover, $a:T\hra G$ induces an embedding $\ov{a}:\ov{T}\hra L_a$ of a maximal elliptic torus of $L_a$.

b) Let $\theta:T(E)\to\B{C}\m$ be a character in general position, trivial on 
$\Ker[T(\C{O})\to\ov{T}(\fq)]$. Denote by $\ov{\theta}:\ov{T}(\fq)\to\B{C}\m$ 
the character of $\ov{T}(\fq)$ defined by $\theta$.
Then there exists a unique irreducible representation $\rho_{a,\theta}$ of $Z(G)(E)G_a$,
whose central character is the restriction of $\theta$, extending the inflation to $G_a$
of the cuspidal Deligne--Lusztig representation $\rho_{\ov{a},\ov{\theta}}$ of $L_a(\fq)$.
We denote by $\pi_{a,\theta}$ the induced cuspidal representation 
$\Ind_{Z(G)(E)G_a}^{G(E)}\rho_{a,\theta}$ of $G(E)$.
\end{Not} 
\begin{Emp} \label{E:coh}
Recall (see \cite[Thm 1.2]{Ko2}) that for every reductive group $G$ over
$E$, $H^1(E,G)$ is canonically isomorphic to the group  $\pi_0(Z(\wh{G})^{\Gm})^D$
of $\B{C}\m$-valued characters of $\pi_0(Z(\wh{G})^{\Gm})$. If $T$ is a maximal torus 
of $G$, we get a commutative diagram:
\[
\CD
H^1(E,T)@>{\sim}>> \pi_0(\wh{T}^{\Gm})^D\\
@VVV          @VVV\\
 H^1(E,G)@>{\sim}>> \pi_0(Z(\wh{G})^{\Gm})^D.
\endCD
\]
In particular, we have a canonical surjection
\[
\wh{T}^{\Gm}/Z(\wh{G})^{\Gm}\to \Coker[\pi_0(Z(\wh{G})^{\Gm})\to\pi_0(\wh{T}^{\Gm})]\isom
(\Ker\,[H^1(E,T)\to H^1(E,G)])^D.
\]
\end{Emp}

\begin{Not} \label{N:char}
a) To every pair $a, a'$ of stably conjugate embeddings $T\hra G$, one associates the class  
$\inv(a',a)\in\Ker\,[H^1(E,T)\to H^1(E,G)]$. This is the class of a cocycle
$c_{\sigma}=g^{-1}\sigma(g)$, where $g\in G(\ov{E})$ is such that $a'=gag^{-1}$
(compare \cite[4.1]{Ko2}).

b) To each  $\ka\in \wh{T}^{\Gm}/Z(\wh{G})^{\Gm}$, an embedding  $a_0:T\hra G$, and 
a character $\theta$ of $T(E)$ as in \rn{repr}, we associate the invariant distribution 
\[
\chi_{a_0,\ka,\theta}:=e(G)\sum_{a}\lan \inv(a,a_0),\ka\ran \chi(\pi_{a,\theta}).
\]
Here $a$ runs over a set of representatives of conjugacy classes of embeddings which are 
stably conjugate to $a_0$, and $\chi(\pi_{a,\theta})$ denotes the character of $\pi_{a,\theta}$.
\end{Not}

\begin{Not} \label{N:endoscopy}
Each pair $(a,\ka)$, where $a:T\hra G$ is an embedding of a maximal torus of 
$G$ and $\ka$ is an element of $\wh{T}^{\Gm}$, gives rise to an isomorphism class
 $\C{E}_{(a,\ka)}$ of an endoscopic datum of $G$.  
Furthermore,  $\C{E}_{(a,\ka)}$ is elliptic if $a(T)$ is an elliptic torus of $G$ 
(see \cite[\S 7]{Ko1} for the definitions of endoscopic data, and compare \cite[II, 4]{La}).

More precisely, each embedding $\eta:\wh{T}\hra\wh{G}$, whose conjugacy class 
corresponds to the stable conjugacy class of $a$, defines an endoscopic datum 
$\C{E}_{(a,\ka,\eta)}=(s,\rho)$, consisting of a semisimple element 
$s=\eta(\ka)$ of $\wh{G}$ and a homomorphism
$\rho:\Gm\overset{\rho_T}{\lra}\Norm_{\Aut\,\wh{G}}(\eta(\wh{T}))_s/\eta(\wh{T})\overset{\rho'}{\lra}
\Out(\wh{G}_s^0)$. Here $\rho_T$ is induced by the $E$-structure of $T$, and $\rho'$ is induced
by the inclusion $\Norm_{\Aut\,\wh{G}}(\eta(\wh{T}))_s\subset \Norm_{\Aut\,\wh{G}}(\wh{G}^0_s)$.
Moreover, the isomorphism class of $\C{E}_{(a,\ka,\eta)}$, denoted by  $\C{E}_{(a,\ka)}$, 
does not depend on $\eta$. 
%
\end{Not}

\begin{Not} \label{N:end}
For each $\gm\in G^{\sr}(E)$ and $\xi\in\wh{G_{\gm}}^{\Gm}$,

(i) put $\C{E}_{(\gm,\xi)}:=\C{E}_{(a_{\gm},\xi)}$, where $a_{\gm}:G_{\gm}\hra G$ is the inclusion map;

(ii) fix an invariant measure $dg_{\gm}$ on $G_{\gm}(E)$, and put 
\[
O_{\gm}(\phi):=\int_{G(E)/G_{\gm}(E)}f(g\gm g^{-1})\frac{dg}{dg_{\gm}}
\]
for each $\phi=fdg\in\C{S}(G(E))$.

(iii) denote by $\ov{\xi}\in \pi_0(\wh{G_{\gm}}^{\Gm}/Z(\wh{G})^{\Gm})$ the class of $\xi$;

(iv) denote by $O^{\ov{\xi}}_{\gm}\in \C{D}(G(E))$ the sum 
$\sum_{\gm'}\lan\inv(\gm',\gm),\ov{\xi}\ran O_{\gm'}$,
taken over a set of representatives of the conjugacy classes stably conjugate to $\gm$, 
where each $dg_{\gm'}$ is compatible with $dg_{\gm}$.
\end{Not}

\begin{Def} \label{D:stab}
Let $\C{E}$ be an endoscopic datum of $G$.

(i) A measure $\phi\in\C{S}(G(E))$ is called {\em $\C{E}$-unstable} if 
$O^{\ov{\xi}}_{\gm}(\phi)=0$ for all pairs $(\gm,\xi)$ as in \rn{end} for which 
$\C{E}_{(\gm,\xi)}$ is isomorphic to $\C{E}$.

(ii) A distribution $F\in\C{D}(G(E))$ is called {\em $\C{E}$-stable} if $F(\phi)=0$ for all 
$\C{E}$-unstable $\phi\in\C{S}(G(E))$.
\end{Def}



\begin{Thm} \label{T:main}
Assume that $p>\dim\,G^{\der}$. Then for each triple  $(a_0,\ka,\theta)$, the distribution  
$\chi_{a_0,\ka,\theta}$ is $\C{E}_{(a_0,\ka)}$-stable.
\end{Thm}

\begin{Not} \label{N: }
For an endoscopic datum $\C{E}=(s,\rho)$, choose a representative $\wt{s}\in\wh{G}^{\ssc}$ of $s$, and
let $Z(\C{E})$ be the set of $z\in Z(\wh{G}^{\ssc})^{\Gm}$ for which there exists $g\in\wh{G}_s$ commuting with 
$\rho:\Gm\to\Out(\wh{G}^0_s)$ such that $g\wt{s}g^{-1}=z\wt{s}$. 
Then $Z(\C{E})$ is a subgroup of $Z(\wh{G}^{\ssc})^{\Gm}$, depending only on the isomorphism class of $\C{E}$.
\end{Not}

\begin{Def} \label{D:adm}
Let $\C{E}$ be an endoscopic datum of $G$. An inner twisting $\varphi:G\to G'$ is called 
{\em $\C{E}$-admissible} if the corresponding class
$\inv(G',G)\in H^1(E,G^{\ad})\cong (Z(\wh{G}^{\ssc})^{\Gm})^D$ is orthogonal to 
$Z(\C{E})\subset Z(\wh{G}^{\ssc})^{\Gm}$.
\end{Def}

\begin{Def} \label{D:def2}
Let $G$ be a reductive group over $E$, $\C{E}=(s,\rho)$ an elliptic endoscopic datum of $G$, 
and $\varphi:G\to G'$ an $\C{E}$-admissible inner twisting. Fix a triple $(a,a';\ka)$, consisting of 
a pair $a:T\hra G$ and  $a':T\hra G'$ of stably conjugate embeddings of maximal tori, and an element 
$\ka\in\wh{T}^{\Gm}$ such that $\C{E}_{(a,\ka)}\cong\C{E}$.

a) Consider $\phi\in\C{S}(G(E))$ and $\phi'\in\C{S}(G'(E))$. They are called 
{\em $(a,a';\ka)$-\\indistinguishable} if they satisfy the following conditions.

(A) For every $\gm\in G^{\sr}(E)$ and  $\xi\in \wh{G_{\gm}}^{\Gm}$ such that 
$\C{E}_{(\gm,\xi)}\cong\C{E}$ and $O^{\ov{\xi}}_{\gm}(\phi)\neq 0$, 

 (i) there exists $\gm'\in G'(E)$ stably conjugate to $\gm$;

(ii) we have, $O^{\ov{\xi}}_{\gm'}(\phi')=\lan \frac{\gm',\gm;\xi}{a',a;\ka}\ran^{-1}
 O^{\ov{\xi}}_{\gm}(\phi)$.\\
Here $\lan \frac{\gm',\gm;\xi}{a',a;\ka}\ran\in\B{C}\m$ is the invariant $\lan \frac{a_{\gm'},a_{\gm};\xi}{a',a;\ka}\ran$
defined in the Appendix for embeddings $a_{\gm}:G_{\gm}\hra G$ and $a_{\gm'}:G_{\gm}\hra G'$  such that
$a_{\gm}(\gm)=\gm$ and $a_{\gm'}(\gm)=\gm'$.

(B) Condition (A) holds if $G$, $a_0$, $\gm$, $\phi$ are interchanged with $G'$, $a'_0$, $\gm'$, $\phi'$. 

b)  The distributions $F\in\C{D}(G(E))$ and $F'\in\C{D}(G'(E))$ are called  
{\em $(a,a';\ka)$-equivalent} if $F(\phi)=F'(\phi')$ for every two  $(a,a';\ka)$-indistinguishable 
measures $\phi$ and $\phi'$.
\end{Def}

\begin{Rem} \label{R:end}
If $\phi$ is $\C{E}_{(a_0,\ka)}$-unstable, then $\phi$ and $\phi'=0$ are 
$(a_0,a'_0;\ka)$-indistinguishable. Therefore every two $(a_0,a'_0;\ka)$-equivalent
distributions $F$ and $F'$ are $\C{E}_{(a_0,\ka)}$-stable.
\end{Rem}

\begin{Main} 
Assume that  $p>\dim\,G^{\der}$.
Let $\varphi:G\to G'$ be an $\C{E}_{(a_0,\ka)}$-admissible inner twisting.
Let $a'_0:T\hra G'$ be an embedding which is  stably conjugate to $a_0$. 
Then the distributions $\chi_{a_0,\ka,\theta}$ on $G(E)$ and  $\chi_{a'_0,\ka,\theta}$ on $G'(E)$ are 
$(a_0,a'_0;\ka)$-equivalent. 
\end{Main}

\begin{Rem}
a) By \rr{end}, \rt{main} follows from the Main Theorem.

b) We believe that a much smaller bound on $p$ would suffice.
\end{Rem}

\section{Basic ingredients of the argument}
\subsection{A generalization of a theorem of Waldspurger} \label{SS:Wa}

Suppose that we are in the situation of \rd{def2}. Then $\varphi$ induces an inner twisting
$\C{G}\to \C{G}'$. As in \rd{def2}, one can define $(a,a';\ka)$-equivalence of 
$F\in\C{D}(\C{G}(E))$ and $F'\in\C{D}(\C{G}'(E))$.

Fix a nontrivial character $\psi: E\to \B{C}\m$, a nondegenerate $G$-invariant pairing 
$\lan\cdot,\cdot\ran$ on $\C{G}$, and $\varphi$-compatible invariant measures on $\C{G}(E)$ and 
$\C{G}'(E)$.
Then $\varphi$ defines a  nondegenerate $G'$-invariant pairing $\lan\cdot,\cdot\ran'$
on $\C{G}'$. These data define Fourier transforms $F\mapsto \C{F}(F)$ on $\C{G}(E)$ and $\C{G}'(E)$.

\begin{Thm} \label{T:Wa}
The distributions $F\in\C{D}(\C{G}(E))$ and $F'\in\C{D}(\C{G}'(E))$ are  $(a,a';\ka)$-equivalent
if and only if $e(G)\C{F}(F)$ and $e(G')\C{F}(F')$ are  $(a,a';\ka)$-equivalent.
\end{Thm}

The proof is a generalization of that of Waldspurger \cite{Wa}
who treated the case $\phi'=0$ (compare also \cite[Thm. 2.7.1]{KP}, where the stable case is treated).

\subsection{Springer hypothesis}

In the notation of \re{finite}, assume that  $\ov{a}(\ov{\C{T}})(\fq)\subset\C{L}(\fq)$
contains an $L$-regular element $\ov{t}$ [and that $p$ is so large that
the logarithm defines an isomorphism $\log:L_{\un}\isom \C{L}_{\nil}$ between unipotent elements
of $L$ and nilpotent elements of $\C{L}$].
Let $\dt_{\ov{t}}$ be the characteristic function of the $\Ad(L(\fq))$-orbit of $\ov{t}$, and let
$\C{F}(\dt_{\ov{t}})$ be its Fourier transform.

We need the following result of \cite{Ka1}.

\begin{Thm} \label{T:Spr}
For every $u\in L_{\un}(\fq)$, we have 
$$
\Tr\,\rho_{\ov{a},\ov{\theta}}(u)=q^{-(\dim\,L-\dim\,\ov{T})/2}\C{F}(\dt_{\ov{t}})(\log(u)).
$$
\end{Thm}

\subsection{Topological Jordan decomposition}
We will call an element $\gm\in G(E)$ {\em compact} if it generates a relatively compact subgroup of $G(E)$.
We will call an element $\gm\in G(E)$ {\em topologically unipotent} if the sequence $\{\gm^{p^n}\}_n$ 
converges to $1$.
Every topologically unipotent element is compact.

The following result is a rather straightforward generalization of \cite[Lem. 2, p. 226]{Ka2}.

\begin{Lem} \label{L:jor}
For every compact element $\gm\in G(E)$ there exists a unique decomposition $\gm=\dt u$ such that
$\dt$ and $u$ commute, $\dt$ is of finite order prime to $p$, and $u$ is topologically unipotent.
In particular, this decomposition is compatible with conjugation and field extensions.
\end{Lem}

\section{A sketch of proof of the Main Theorem}

\subsection{Reformulation of the problem} 

\begin{Not}
To each $a:T\hra G$ and $\theta:T(E)\to\B{C}\m$ as in \rn{repr} we associate a function $t_{a,\theta}$ on $G(E)$ 
supported on  $Z(G)(E) G_a$ and equal to $\Tr\,\rho_{a,\theta}$ there.
Since $t_{a,\theta}$ is cuspidal the integral
\[
F_{a,\theta}(\gm):=\frac{1}{\mu((G^{\ad})_a)} 
\int_{G(E)/Z(G)(E)}t_{a,\theta}(g\gm g^{-1})dg
\] 
stabilizes for every $\gm\in G^{\sr}(E)$ (see \cite[Lem. 23]{HC}), thus providing us with a locally constant 
invariant function $F_{a,\theta}$ on $G^{\sr}(E)$.
\end{Not}

\begin{Lem} \label{L:HCH}
For each $a$ and $\theta$, $F_{a,\theta}$ is a locally $L^1$-function on $G(E)$.
Moreover, the corresponding distribution equals $\chi({\pi_{a,\theta}})$.
\end{Lem} 
\begin{proof}
The assertion follows from Harish-Chandra's theorem \cite[Thm. 16]{HC}.
\end{proof}

\begin{Not} \label{N:sigma}
For every $\gm_0\in G^{\sr}$ and 
$\ov{\xi}\in\pi_0(\wh{G_{\gm_0}}^{\Gm}/Z(\wh{G})^{\Gm})$ we define
\begin{equation} \label{E:redst}
\Sigma_{G;\gm_0,\ov{\xi};a_0,\ka}:=e(G)\sum_{a}\sum_{\gm} 
\lan \inv(a,a_0),\ka\ran \lan \inv(\gm,\gm_0),\ov{\xi}\ran ^{-1}F_{a,\theta}(\gm),
\end{equation}
where $a$ and $\gm$ run over sets of representatives of the conjugacy classes within the 
stable conjugacy classes of $a_0$ and $\gm_0$, respectively. 
\end{Not}

\begin{Thm} \label{T:main'}
For all $\gm_0\in G^{\sr}$ and $\ov{\xi}\in\pi_0(\wh{G_{\gm_0}}^{\Gm}/Z(\wh{G})^{\Gm})$ 
such that $\Sigma_{G;\gm_0,\ov{\xi};a_0,\ka}\neq 0$, 

 $(i)$ there exists a representative $\xi\in \wh{G_{\gm_0}}^{\Gm}$ of $\ov{\xi}$ such that
$\C{E}_{(\gm_0,\xi)}\cong\C{E}_{(a_0,\ka)};$

 $(ii)$ if $\varphi:G\to G'$ is $(\C{E},a_{\gm_0},\ov{\xi})$-admissible (see \rd{app}), 
then for every $\xi$ as in  $(i)$ and every stably conjugate $\gm'_0\in G'(E)$ of $\gm_0$ 
we have \[\Sigma_{G';\gm'_0,\ov{\xi};a'_0,\ka}=
\lan \frac{\gm'_0,\gm_0;\xi}{a',a;\ka}\ran\Sigma_{G;\gm_0,\ov{\xi};a_0,\ka}.\]
\end{Thm}

\begin{Emp}
It follows from \rl{HCH} that \rt{main'} is equivalent to the Main Theorem.
Moreover, by standard arguments, \rt{main'} reduces to the case when the derived group of 
$G$ is simply connected. 
\end{Emp}

\begin{Emp} \label{E:red}
From now on we will assume that $G^{\der}=G^{\ssc}$. 
In particular, the centralizer of each semisimple element of $G$ is 
connected, and each $G_a$ is a maximal compact subgroup of $G(E)$.
We fix $(\gm_0,\ov{\xi})$ such that $\Sigma_{G;\gm_0,\ov{\xi};a_0,\ka}\neq 0$.
Since $\Sigma_{G;z\gm_0,\ov{\xi};a_0,\ka}=\theta(z)\Sigma_{G;\gm_0,\ov{\xi};a_0,\ka}$
for each  $z\in Z(G)(E)$ and since the support of each $t_{a,\theta}$ consists of elements compact modulo center,
we can assume that $\gm_0$ is compact with topological Jordan decomposition  $\gm_0=\dt_0 u_0$.
Moreover, we can assume either that $\gm_0$ is topologically unipotent, or that $\dt_0\notin Z(G)(E)$.
\end{Emp}

\subsection{The topologically unipotent case} \label{SS:topun} 
\begin{Emp}
Since $p$ does not divide the order of $Z(G^{\der})$, the canonical map
$G^{\der}(E)_{\tu}\times Z(G)(E)_{\tu}\to G(E)_{\tu}$ is an isomorphism.
Therefore to prove \rt{main'} for topologically unipotent $\gm_0$, we can assume that 
$G$ is semisimple and simply connected.
\end{Emp}

\begin{Not}
Denote by $\Phi_{G}:G\to\C{G}$ the composition map
\[
G\overset{\Ad}{\lra} GL(\C{G})\overset{\log_{(p)}}{\lra}\End(\C{G})\overset{\pr}{\lra}\C{G},
\]
where $\log_{(p)}(1-A)=-\sum_{i=1}^{p-1} \frac{A^i}{i}$,
and $\pr$ is the canonical projection, defined by the standard pairing 
$(A,B)\mapsto\Tr AB$ on $\End(\C{G})$.
\end{Not}

\begin{Lem} \label{L:bij}
The map $\Phi_G$ defines a $G(E)$-equivariant homeomorphism \\
$G(E)_{\tu}\isom \C{G}(E)_{\tn}$, where $G(E)$ acts by conjugation
Moreover, for every parahoric subgroup $G_x$ of $G(E)$, 
$\Phi_G$ induces a bijection $(\Phi_G)_x:(G_x)_{\tu}\isom (\C{G}_x)_{\tn}$, which in turn
induces the logarithm map $\log:(L_x)_{\un}(\fq)\isom (\C{L}_x)_{\nil}(\fq)$.
\end{Lem}

\begin{Not} 
a) By our assumption on $p$, there exists $t\in\C{T}(\C{O})$ whose reduction 
$\ov{t}\in \ov{\C{T}}(\fq)$ is not fixed by any nontrivial element of Weyl group of $G$.

b) For every $a:T\hra G$ as in \rn{repr}, we denote by $\Om_{a,t}\subset \C{L}_a(\fq)$ 
the $\Ad(L_a(\fq))$-orbit of $\ov{a}(\ov{t})$, by $\wt{\Om}_{a,t}\subset \C{G}_{a}\subset \C{G}$
the preimage of $\Om_{a,t}$, and let $\dt_{a,t}$ be the characteristic function of $\wt{\Om}_{a,t}$.

c) As  the centralizer $G_y$ of each $y\in\wt{\Om}_{a,t}$ is $G_a$-conjugate to $a(T)$, the integral
\[
\Dt_{a,t}(x):=\frac{1}{\mu(G_a)}\int_{G(E)}\dt_{a,t}(\Ad(g)x) dg
\]
converges absolutely for each $x\in\C{G}(E)$. Thus it defines an element of $\C{D}(\C{G}(E))$. 
Similarly to \rn{char}, we consider  
${\Dt}_{a_0,\ka,t}:=e(G)\sum_{a}\lan \inv(a,a_0),\ka\ran \Dt_{a,t}\in \C{D}(\C{G}(E))$.
\end{Not}

\begin{Lem} \label{L:four}
Let $\C{I}^+\subset\C{G}(E)$ be a maximal topologically nilpotent subalgebra.
Assume that $\psi:E\to\B{C}\m$ is trivial on the maximal ideal $M\subset\C{O}$ and
induces a nontrivial character of $\fq$. Then for each $u\in G(E)_{\tu}$ we have
\[t_{a,\theta}(u)=\mu(\C{I}^+)^{-1} \C{F}(\dt_{a,t})(\Phi_G(u)).\]
\end{Lem}
\begin{proof}
The assumption on $\psi$ implies that $\C{G}_{a^+}$ is the orthogonal complement of
$\C{G}_{a}$ with respect to the pairing $(x,y)\mapsto\psi(\lan x,y\ran)$. Therefore our lemma
is an immediate consequence of the definition of the Fourier transform (over $E$ and $\fq$),
 \rt{Spr}, \rl{bij}, and the equality $q^{(\dim\,L_a-\dim\,\ov{T})/2}\mu(\C{G}_{a^+})=\mu(\C{I}^+)$.
\end{proof}

\begin{Emp}
Now we are ready to show that $(\chi_{a_0,\ka,\theta})_{|G(E)_{\tu}}$ and 
$(\chi_{a'_0,\ka,\theta})_{|G'(E)_{\tu}}$ are $(a_0,a'_0;\ka)$-equivalent. First of all, by 
direct calculations, $e(G){\Dt}_{a_0,\ka,t}$ is $(a_0,a'_0;\ka)$-equivalent to $e(G'){\Dt}_{a'_0,\ka,t}$. 
Hence by \rt{Wa},  $\C{F}({\Dt}_{a_0,\ka,t})$ is $(a_0,a'_0;\ka)$-equivalent to 
$\C{F}({\Dt}_{a'_0,\ka,t})$. Using Lemmas \ref{L:HCH}, \ref{L:bij} and \ref{L:four} we see that 
$\chi_{a_0,\ka,\theta}$ has the same restriction to $G(E)_{\tu}$ as 
$\mu(\C{I}^+)^{-1}\Phi_G^*(\C{F}({\Dt}_{a_0,\ka,t}))$, and similarly for $\chi_{a'_0,\ka,\theta}$ and 
$\mu(\C{I'}^+)^{-1}\Phi_{G'}^*(\C{F}({\Dt}_{a'_0,\ka,t}))$.
Since $\Phi_G$ is $G^{\ad}$-invariant algebraic morphism defined over
$E$, the assertion follows from the equality $\mu(\C{I}^+)=\mu(\C{I'}^+)$. 
\end{Emp}

\subsection{The general case}

It remains to prove \rt{main'} for $\dt_0\notin Z(G)(E)$ (see \re{red}).
We are going to deduce the assertion from that for $G_{\dt_0}$.

\begin{Prop} \label{P:reduction}
For every embedding $a:T\hra G$ and a compact element $\gm$ in $G(E)$ with topological Jordan decomposition
$\gm=\dt u$, we have 
\[
e(G)F_{a,\theta}(\gm)=e(G_{\dt})\sum_{b}\theta(b^{-1}({\dt}))F_{b,\theta}(u).
\]
Here $b$ runs over the set of conjugacy classes of embeddings $b:T\hra G_{\dt}$
whose composition with the inclusion $G_{\dt}\subset G$ is conjugate to $a$. 
\end{Prop}
\begin{proof}
The proposition follows by direct calculation from the recursive formula
(\cite[Thm. 4.2]{DL}) for characters of Deligne--Lusztig representations.
\end{proof}

\begin{Not}
a) We say that $t\in T(E)$ is {\em $(G,a_0,\gm_0)$-relevant} if there exists an embedding 
$b_0:T\hra  G_{\dt_0}\subset G$ stably conjugate to $a_0$ such that $b(t)=\dt_0$. 

b) Assume that  $t\in T(E)$ is $(G,a_0,\gm_0)$-relevant.
Since $a_0(T)\subset G$ is elliptic, for each $\dt\in G(E)$ stably conjugate to $\dt_0$ 
there exists an embedding $b_{t,\dt}:T\hra G_{\dt}\subset G$ stably conjugate to $a_0$ such that 
$b_{t,\dt}(t)=\dt$. Further, $b_{t,\dt}$ is unique up to stable conjugacy, and 
the endoscopic datum $\C{E}_{t,\ka}:=\C{E}_{(b_{{\dt},t},\ka)}$ of $G_{\dt_0}$ is independent of $\dt$.

c) We will write $\dt_1\underset{\C{E}_{t,\ka}}\sim\dt$ (resp. $\dt'\underset{\C{E}_{t,\ka}}\sim\dt$) if
$\dt,\dt_1\in G(E)$ (resp. $\dt\in G(E)$ and $\dt'\in G'(E)$) are stably conjugate to $\dt_0$, 
and $G_{\dt_1}$ (resp.  $G'_{\dt'}$) is an  $\C{E}_{t,\ka}$-admissible inner form of $G_{\dt}$
 (see \rd{adm}).
\end{Not}

\begin{Emp}
Using \rp{reduction}, we see that 
\begin{equation} \label{E:for2}
\Sigma_{G;\gm_0,\ov{\xi};a_0,\ka}=\sum_{t}\theta(t)\sum_{{\dt}}I_{t,\dt},
\end{equation}
where

 (i) $t$ runs over the set of $(G,a_0,\gm_0)$-relevant elements of $T(E)$;

 (ii) ${\dt}$ runs over a set of representatives of the conjugacy classes within the stable conjugacy class
of ${\dt}_0$;

(iii) $I_{t,\dt}$ vanishes unless there exists an element  $\gm\in G(E)$ stably conjugate to $\gm_0$
with topological Jordan decomposition $\gm={\dt}u$, in which case we get 
\begin{equation} \label{E:for4}
I_{t,\dt}=\lan \inv(b_{t,\dt},a_0),\ka\ran \lan \inv(\gm,\gm_0),\ov{\xi}\ran ^{-1}
\Si_{G_{\dt};u,\ov{\xi};b_{t,\dt},\ka}.
\end{equation}
\end{Emp}

\begin{Emp}
For simplicity of the exposition, we will restrict ourselves to the case when 
$\gm_0\in G(E)$ is elliptic. 
Choose $t$ which has a nonzero contribution to (\ref{E:for2}). 
Replacing $\dt_0$ by a stably conjugate element we can assume that
$\sum_{\dt\underset{\C{E}_{t,\ka}}\sim\dt_0}I_{t,\dt}\neq 0$ and $I_{t,\dt_0}\neq 0$. 
So $\Si_{G_{\dt_0};u_0,\ov{\xi};b_{t,\dt},\ka}\neq 0$.
Hence by \rt{main'} for $G_{\dt_0}$ there exists a representative $\xi\in \wh{G_{\gm_0}}^{\Gm}$ of 
$\ov{\xi}$ such that the endoscopic datum $\C{E}_{(u_0,\xi)}$ of $G_{\dt_0}$ is isomorphic to $\C{E}_{t,\ka}$.  
Therefore there exist embeddings $\eta_1:\wh{G_{\gm_0}}\hra\wh{G_{\dt_0}}$ and 
$\eta_2:\wh{T}\hra\wh{G_{\dt_0}}$ such that 
$\C{E}_{(\gm_0,\xi,\eta_1)}=\C{E}_{(b_{t,\dt_0},\ka,\eta_2)}$ (compare \rn{endoscopy}) and 
$\eta_2(\ka)=z\eta_1(\xi)$ for a certain $z\in Z(\wh{G_{\dt_0}})^{\Gm}$. 
Moreover, $z$ is defined up to multiplication by an element of $Z(\C{E}_{t,\ka})$. Therefore 
for all $\dt\underset{\C{E}_{t,\ka}}\sim\dt_0$, the expression $\lan \inv(\dt,\dt_0),z\ran$ is 
independent of the choice of the $\eta_i$'s.
\end{Emp}

\begin{Cl} \label{C:eq}
For each $\dt\underset{\C{E}_{t,\ka}}\sim\dt_0$ we have 
$I_{t,\dt}=\lan \inv(\dt,\dt_0),z\ran I_{t,\dt_0}$.
\end{Cl}
\begin{proof}
Since $\Si_{G_{\dt_0};u_0,\ov{\xi};b_{t,\dt},\ka}\neq 0$, \rt{main'} for inner forms $G_{\dt}$ and $G_{\dt_0}$
implies that for every stably conjugate $u\in G_{\dt}(E)$ of $u_0\in G_{\dt_0}(E)$, we have
\[
\Sigma_{G_{\dt};u,\ov{\xi};b_{t,\dt},\ka}=\lan \frac{u,u_0;\xi}{b_{t,\dt} ,b_{t,\dt_0};\ka}\ran
\Sigma_{G_{\dt_0};u_0,\ov{\xi};b_{t,\dt_0},\ka}.
\]
Then $\gm:=\dt u\in G(E)$ is stably conjugate to $\gm_0$, and the assertion follows by 
direct calculations from (\ref{E:for4}).
\end{proof}

\begin{Emp}
Now we are ready to show the validity of (i),(ii) of \rt{main'}.

(i) As $\sum_{\dt\underset{\C{E}_{t,\ka}}\sim\dt_0}I_{t,\dt}\neq 0$, we get from \rcl{eq} that
$\sum_{\dt\underset{\C{E}_{t,\ka}}\sim\dt_0}\lan \inv(\dt,\dt_0),z\ran\neq 0$. By the definition of 
$\C{E}_{t,\ka}$-equivalence, this implies that $z$ belongs to $Z(\C{E}_{t,\ka})Z(\wh{G})^{\Gm}$.
Thus changing $\eta_1$ (or $\eta_2$), we can assume that $z\in Z(\wh{G})^{\Gm}$.
Since  $(\gm_0,\xi)$ and $(b_{t,\dt_0},\ka)$ define isomorphic endoscopic data of $G_{\dt_0}$, 
we therefore conclude that $\C{E}_{(\gm_0,\xi)}\cong\C{E}_{(a_0,\ka)}$.


(ii) Since $T$ is elliptic, an element $t\in T(E)$ is $(G,a_0,\gm_0)$-relevant if and only if it is 
 $(G',a'_0,\gm'_0)$-relevant. Thus it will suffice to show that for each such $t$ we have
$\sum_{\dt'}I_{t,\dt'}=\lan\frac{\gm'_0,\gm_0;\xi}{a'_0,a_0;\ka}\ran \sum_{\dt}I_{t,\dt}$.
For every stably conjugate $\dt\in G(E)$ of $\dt_0$, there exists
a stably conjugate $\dt'\in G'(E)$ of $\dt'_0$ such that  $\dt'\underset{\C{E}_{t,\ka}}\sim\dt$.
Therefore it will suffice to show that for every such pair $\dt'\underset{\C{E}_{t,\ka}}\sim\dt$, we have 
$I_{t,\dt'}=\lan\frac{\gm'_0,\gm_0;\xi}{a'_0,a_0;\ka}\ran I_{t,\dt}$.
The latter equality can be proved by the same arguments as \rcl{eq}.
\end{Emp}

\appendix
\section{}

Let $G$ be a reductive group over $E$, $\C{E}=(s,\rho)$ an elliptic endoscopic datum of $G$,
and $\varphi:G\to G'$ an $\C{E}$-admissible inner twisting. To every two triples 
$(a'_i,a_i;\ka_i)$, $i=1,2$, where $a_i:T_i\hra G$ and  $a'_i:T_i\hra G'$ 
are stably conjugate embeddings of maximal tori, and $\ka_i$ is an element $\wh{T_i}^{\Gm}$
such that $\C{E}_{(a_i,\ka_i)}$ is isomorphic to $\C{E}$, we are going to define an invariant 
$\lan \frac{a'_1,a_1;\ka_1}{a'_2,a_2;\ka_2}\ran\in\B{C}\m$.

{\bf Step 1.} Replacing $G$, $G'$, $T_i$, $\ka_i$ and $\C{E}$ by $G^{\ssc}$, $G'^{\ssc}$,  
$T^{\ssc}_i:=a_i^{-1}(G^{\ssc})=a'^{-1}_i(G'^{\ssc})$, the image of $\ka_i$ in  
$\wh{T^{\ssc}_i}^{\Gm}$, and the corresponding endoscopic datum of 
$G^{\ssc}$ respectively, we can assume that $G$ is semisimple and simply connected.
Let $T_{1,2}$ be the quotient 
$[T_1\times T_2]/Z(G)$.

{\bf Step 2.}  Choose elements $g_1,g_2$ and $\{\wt{c}_{\si}\}_{\si\in\Gm}$ of $G(\ov{E})$ such that 
$a'_i=\varphi(g_i a_i g_i^{-1})$ and each $\wt{c}_{\si}$ is a representative of 
$\varphi^{-1}{}^{\si}\varphi\in G^{\ad}(\ov{E})$. Then  each 
$g_i^{-1}\wt{c}_{\si}{}^{\si}g_i\in G(\ov{E})$ belongs to $a_i(T_i(\ov{E}))$, and the images of 
$(a_1^{-1}(g_1^{-1}\wt{c}_{\si}{}^{\si}g_1),
a_2^{-1}(g_2^{-1}\wt{c}_{\si}{}^{\si}g_2))$ in 
$T_{1,2}(\ov{E})$ form a cocycle, whose cohomology class  
$\inv(a'_1,a_1;a'_2,a_2)\in  H^1(E,T_{1,2})$ is independent of the choices.

{\bf Step 3.} Choose embeddings $\eta_i:\wh{T_i}\hra\wh{G}$ such that 
$\C{E}_{(a_i,\ka_i,\eta_i)}=(s,\rho)$ and a representative $\wt{s}\in\wh{G}^{\ssc}=\wh{G^{\ad}}$ of $s$. 
Put $T^{\ad}_i:=T_i/a_i^{-1}(Z(G))$. Each $\eta_i$ defines 
an embedding $\wt{\eta}_i:\wh{T^{\ad}_i}\hra\wh{G^{\ad}}$, hence an 
element $\wt{\ka}_i=\ka(\wt{s},\eta_i):=\wt{\eta}_i^{-1}(\wt{s})\in\wh{T^{\ad}_i}$. 
Then the image of $(\wt{\ka}_1,\wt{\ka}^{-1}_2)$ in 
$\wh{T^{\ad}_1}\times \wh{T^{\ad}_2}/Z(\wh{G^{\ad}})\cong\wh{T_{1,2}}$, 
denoted by $\frac{\ka_1}{\ka_2}$, is $\Gm$-invariant. 
Moreover, as $\varphi:G\to G'$ is $\C{E}$-admissible, the expression 
$\lan \frac{a'_1,a_1;\ka_1}{a'_2,a_2;\ka_2}\ran:=
\lan\inv(a'_1,a_1;a'_2,a_2),\frac{\ka_1}{\ka_2}\ran\in\B{C}\m$ 
is independent of the choices.

\begin{Def} \label{D:app}
Let $\C{E}=(s,\rho)$ be an endoscopic datum of $G$, 
$\varphi:G\to G'$ an inner twisting, $a:T\hra G$ an embedding 
of a maximal torus, and $\ka$ an element of $\wh{T}^{\Gm}$ such that
$\C{E}_{(a,\ka)}\cong\C{E}$.
We say that $\varphi:G\to G'$ is  
{\em $(\C{E},a,\ov{\ka})$-admissible}, if for all representatives 
$\ka'\in\wh{T}^{\Gm}$ of $\ov{\ka}\in\pi_0(\wh{T}^{\Gm}/Z(\wh{G})^{\Gm})$ 
satisfying $\C{E}_{(a,\ka')}\cong\C{E}$, all embeddings 
$\eta,\eta':\wt{T}\hra\wh{G}$ such that 
$\C{E}_{(a,\ka,\eta)}=\C{E}_{(a,\ka',\eta')}=(s,\rho)$ and all representatives
 $\wt{s}\in\wh{G^{\ad}}$ of $s$,  
the difference $\ka(\wt{s},\eta')-\ka(\wt{s},\eta)\in Z(\wh{G^{\ad}})^{\Gm}$ 
is orthogonal to
$\inv(G',G)\in H^1(E,G^{\ad})$.
\end{Def}

\begin{Rem} \label{R:app} 
a) Every $(\C{E},a,\ov{\ka})$-admissible inner twisting is $\C{E}$-admissible.

b) If $a(T)\subset G$ is elliptic, then every $\C{E}$-admissible 
inner twisting  is $(\C{E},a,\ov{\ka})$-admissible.

c) $\varphi:G\to G'$ is $(\C{E},a_1,\ov{\ka}_1)$-admissible if and only if 
$\lan \frac{a'_1,a_1;\ka_1}{a'_2,a_2;\ka_2}\ran=
\lan\frac{a'_1,a_1;\ka'_1}{a'_2,a_2;\ka_2}\ran$  for all  representatives
$\ka'_1\in\wh{T}_1^{\Gm}$ of $\ov{\ka}_1\in\pi_0(\wh{T}_1^{\Gm}/Z(\wh{G})^{\Gm})$ 
satisfying $\C{E}_{(a_1,\ka'_1)}\cong\C{E}$. 
\end{Rem}

\end{document}